# PROJECTIVITY AND ISOMORPHISM OF STRICTLY SIMPLE ALGEBRAS

KEITH A. KEARNES AND ÁGNES SZENDREI

ABSTRACT. We describe a sufficient condition for the localization functor to be a categorical equivalence. Using this result we explain how to simplify the test for projectivity. This leads to a description of the strictly simple algebras which are projective in the variety they generate. A byproduct of our efforts is the result that if **A** and **B** are strictly simple and generate the same variety, then $\mathbf{A} \cong \mathbf{B}$ or else both are strongly abelian.

## 1. INTRODUCTION

Let **A** be an algebra and let $e$ be a unary term in the language of **A**. The term $e$ is **idempotent** if $\mathbf{A} \models e^2 = e$. If $e$ is idempotent, we call the set $e(A)$ a **neighborhood** of **A**. In Section 2 we explain how to localize the structure of **A** to the neighborhood $e(A)$. If $\mathcal{V}$ is a variety of algebras and $\mathcal{V} \models e^2 = e$, then localization to the range of $e$ is a functorial construction on the members of $\mathcal{V}$.

In this paper we analyze properties of the localization functor. In Section 2 we describe broad sufficient conditions which guarantee that this functor is a categorical equivalence. We use this in the succeeding sections to reduce questions about projectivity and isomorphism to the case of term minimal algebras. As an application, we show that if **A** and **B** are strictly simple algebras which generate the same variety and are not strongly abelian then $\mathbf{A} \cong \mathbf{B}$. It follows then that if **A** and **B** are finite simple algebras of type **2** which generate the same variety then $\mathbf{A} \cong \mathbf{B}$. This settles a question left open in [4] where it is demonstrated that the analogous statement holds true for types **3** and **4** and fails for types **1** and **5**.

## 2. THE LOCALIZATION FUNCTOR

Let **A** be an algebra and let $e$ be an idempotent of **A**. For every term $t$ of **A**, $et$ is a term of **A** of the same arity as $t$ such that the neighborhood $e(A)$ is closed under the term operation $et$. We will use the symbol $e(\mathbf{A})$ to denote the following algebra. The universe of $e(\mathbf{A})$ is the neighborhood $e(A)$. The set of fundamental operation symbols will be the set $\{et \colon t \text{ is a term in the language of } \mathbf{A}\}$, and the interpretation of $et$ as an operation on $e(A)$ is the obvious one: the restriction $et|_{e(A)}$ of the term operation $et$ of $\mathbf{A}$ to $e(A)$. The algebra $e(\mathbf{A})$ is called the **localization of A** to the neighborhood $e(A)$.

Fixing the similarity type of $e(\mathbf{A})$ as we have allows us to consider the localization construction $\mathbf{A} \mapsto e(\mathbf{A})$ for any class $\mathcal{K}$ of similar algebras in which the identity $e^2 = e$ holds. As a result, we get a class $e(\mathcal{K}) = \{e(\mathbf{A}) \colon \mathbf{A} \in \mathcal{K}\}$ of similar algebras.

There is a natural way to extend the object mapping
$$e \colon \mathcal{K} \to e(\mathcal{K}), \quad \mathbf{A} \mapsto e(\mathbf{A})$$
to a functor: to each homomorphism $\varphi \colon \mathbf{B} \to \mathbf{C}$ we define the corresponding homomorphism to be
$$e(\varphi) = \varphi|_{e(B)} \colon e(\mathbf{B}) \to e(\mathbf{C}).$$
Here $e(\varphi)$ is not only a handy notation for the image of $\varphi$ under the functor $e$; it can also be interpreted as the image of the subalgebra $\varphi$ of $\mathbf{B} \times \mathbf{C}$ under the term operation $e$. It is straightforward to check that with this latter interpretation we have $e(\varphi) = \varphi|_{e(B)}$, that $\varphi|_{e(B)}$ is indeed a homomorphism from $e(\mathbf{B})$ into $e(\mathbf{C})$ and that $e \colon \mathcal{K} \to e(\mathcal{K})$ is a functor. It is not hard to show that if $\mathcal{K}$ is closed under the formation of isomorphic images, subalgebras, products or ultraproducts, then so is $e(\mathcal{K})$. If $\mathcal{K}$ is closed under the formation of homomorphic images of subalgebras, then so is $e(\mathcal{K})$. Thus, if $\mathcal{K}$ is a variety, quasivariety, prevariety or pseudovariety, then so is $e(\mathcal{K})$.

There are two concepts about the relationship between an idempotent $e$ and an algebra $\mathbf{A}$ which we shall find interesting. To define the first concept, let $\mathbf{A}$ be any algebra and let $e$ be any idempotent unary term of $\mathbf{A}$. We say that $e$ **separates A**, or $e$ **is separating for A**, provided that for every $a \neq b$ in $\mathbf{A}$ there is a unary term $g$ such that $eg(a) \neq eg(b)$. Separation will be a basic concept in this paper, so let us prove a simple characterization of what it means for $e$ to separate $\mathbf{A}$.

**Lemma 2.1.** *Let $\mathbf{A}$ be an algebra and let $e$ be an idempotent term of $\mathbf{A}$. Then $e$ separates $\mathbf{A}$ if and only if*
  (1) *any isomorphism between subalgebras of $\mathbf{A}$ which restricts to the identity on $e(A)$ is the identity on its domain, and*
  (2) *any congruence on a subalgebra of $\mathbf{A}$ which restricts trivially to $e(A)$ is trivial.*

*Proof.* For the purposes of this proof, call a pair of elements $a, b \in A$ an **inseparable pair** if $eg(a) = eg(b)$ holds for every unary term $g$.

Clearly $e$ separates $\mathbf{A}$ if and only if $\mathbf{A}$ has no inseparable pair of distinct elements.

Assume that $e$ separates $\mathbf{A}$, and therefore that $\mathbf{A}$ has no inseparable pair of distinct elements. If $\varphi\colon \mathbf{B} \to \mathbf{C}$ is an isomorphism between subalgebras of $\mathbf{A}$ which is the identity on $e(A)$ then $(b, \varphi(b))$ is an inseparable pair for all $b \in B$. Since we have assumed that $e$ separates $\mathbf{A}$ we get that $\varphi = \mathrm{id}_{\mathbf{B}}$, so (1) holds. If $\alpha$ is a congruence on a subalgebra whose restriction to $e(A)$ is trivial, then any $(b, c) \in \alpha$ is an inseparable pair. The fact that $e$ separates $\mathbf{A}$ implies that $\alpha$ is trivial, so (2) holds.

Now assume that (1) and (2) hold. Let $\theta$ denote the equivalence relation on $\mathbf{A}$ comprised of the inseparable pairs. Observe that for every $(a, b) \in \theta$ we have $e(a) = e(b)$, so $e(\theta)$ is the equality relation on $e(A)$. Choose $(c, d) \in \theta$ and let $S$ denote the subuniverse of $\mathbf{A}^2$ generated by $(c, d)$. A typical member of $S$ is of the form $(g(c), g(d))$ where $g$ is a unary term so, since $(c, d) \in \theta$, we get that $S \subseteq \theta$ as well. From this we get that $R_1 := S^{\cup} \circ S$ and $R_2 := S \circ S^{\cup}$ are symmetric subuniverses of $\mathbf{A}^2$ contained in $\theta$ which have the property that

$$(a, b) \in R_i \implies (a, a), (b, b) \in R_i.$$

The transitive closure of each $R_i$ is a congruence on a subalgebra of $\mathbf{A}$ which, since $R_i \subseteq \theta$, restricts trivially to $e(A)$. We have assumed that (2) holds, so each $R_i$ is the equality relation on a subalgebra. But the statement that $S^{\cup} \circ S$ and $S \circ S^{\cup}$ are subalgebras of $\mathbf{A}^2$ which are subsets of the equality relation on $A$ means precisely that $S$ is the graph of an isomorphism between subalgebras of $\mathbf{A}$. Since $S \subseteq \theta$ this isomorphism is the identity on $e(A)$, so $S$ is the identity relation on its domain. This is true for $S$ generated by an arbitrarily chosen $(c, d) \in \theta$, so $\theta$ is the equality relation on $\mathbf{A}$. Thus, $\mathbf{A}$ contains no inseparable pair of distinct elements, hence $e$ separates $\mathbf{A}$. $\square$

Now we turn to the second concept that interests us. As usual, let $\mathbf{A}$ be an algebra and let $e$ be an idempotent unary term of $\mathbf{A}$. We say that $e$ is **dense** for $\mathbf{A}$ if $\mathbf{A}$ is generated as an algebra by $e(A)$.

A basic lemma relating the properties of separation and density of the term $e$ to the properties of the functor $e$ is the following.

**Lemma 2.2.** *Let $\mathcal{V}$ be a variety and let $e$ be a unary term in the language of $\mathcal{V}$ for which $\mathcal{V} \models e^2 = e$.*

(1) *The subclass of all $\mathbf{A} \in \mathcal{V}$ for which $e$ is separating is a prevariety. (That is, it is a class closed under the formation of isomorphic images, subalgebras and products.)*

(2) *Consider $e$ as a functor from $\mathcal{V}$ to $e(\mathcal{V})$. Choose $\mathbf{A}, \mathbf{B} \in \mathcal{V}$ and consider the induced mapping $\mathrm{Hom}(\mathbf{A}, \mathbf{B}) \longrightarrow \mathrm{Hom}(e(\mathbf{A}), e(\mathbf{B}))$.*
  (i) *If $e$ is dense for $\mathbf{A}$, then this mapping is injective.*
  (ii) *If $e$ is dense for $\mathbf{A}$ and separating for $\mathbf{B}$, then this mapping is surjective.*

*Proof.* The first claim of the lemma is proved by noting that $e$ is separating for $\mathbf{A} \in \mathcal{V}$ if and only if

$$\mathbf{A} \models \bigwedge_g eg(x) = eg(y) \Longrightarrow x = y.$$

Implications of this form, which may have infinitely many conjuncts, are preserved by the formation of isomorphic images, subalgebras and products. Therefore, the collection of members of $\mathcal{V}$ for which $e$ is separating is a prevariety.

Next we prove statement $(2)(i)$. Assume that $\varphi, \psi \in \mathrm{Hom}(\mathbf{A}, \mathbf{B})$ and that $e(\varphi) = e(\psi)$. This means exactly that $\varphi$ and $\psi$ have the same restriction to $e(A)$. The equalizer of $\varphi$ and $\psi$ is a subalgebra of $\mathbf{A}$ which, since $e(\varphi) = e(\psi)$, contains $e(A)$. Since $e$ is dense for $\mathbf{A}$ we get that the equalizer of $\varphi$ and $\psi$ is $\mathbf{A}$, and so $\varphi = \psi$.

Finally we must prove that $(2)(ii)$ holds. For this, choose $\lambda \in \mathrm{Hom}(e(\mathbf{A}), e(\mathbf{B}))$. We need to prove that there is a $\widehat{\lambda} \in \mathrm{Hom}(\mathbf{A}, \mathbf{B})$ such that $e(\widehat{\lambda}) = \lambda$. We identify a homomorphism with its graph, thus $y = \lambda(x)$ is synonymous with $(x, y) \in \lambda$. Stipulating this, $\lambda$ is a subuniverse of $e(\mathbf{A}) \times e(\mathbf{B})$. Let $\widehat{\lambda}$ be the subuniverse of $\mathbf{A} \times \mathbf{B}$ generated by $\lambda$.

We claim that $e(\widehat{\lambda}) = \lambda$. This can be justified as follows. From the definition we have $\lambda \subseteq \widehat{\lambda}$, so $\lambda = e(\lambda) \subseteq e(\widehat{\lambda})$. Conversely,

$$\begin{aligned}(a, b) \in e(\widehat{\lambda}) &\Longrightarrow (a, b) \in \{et((r_0, s_0), \dots) \mid (r_0, s_0), \dots \in \lambda\} \\ &\Longrightarrow (a, b) \in \lambda,\end{aligned}$$

so $e(\widehat{\lambda}) = \lambda$. (The first implication follows from the definition of $\widehat{\lambda}$ and the second implication follows from the fact that $et$ is a term in the language of $e(\mathcal{V})$ and that $\lambda$ is a subuniverse of $e(\mathbf{A}) \times e(\mathbf{B})$.) Therefore, if we prove that $\widehat{\lambda}$ is a mapping from $A$ to $B$ then we will have finished the proof of $(2)(ii)$.

Assume that $(a, b), (a, c) \in \widehat{\lambda} \subseteq \mathbf{A} \times \mathbf{B}$ and that $b \neq c$. Since $e$ separates $\mathbf{B}$ there is a unary term $g$ such that $eg(b) \neq eg(c)$. But now we have that $(eg(a), eg(b)), (eg(a), eg(c)) \in \lambda$ and $eg(b) \neq eg(c)$, which contradicts the fact that $\lambda$ is a mapping. We conclude that if $(a, b), (a, c) \in \widehat{\lambda}$ then $b = c$. This implies that $\widehat{\lambda}$ is a partial homomorphism from $\mathbf{A}$ to $\mathbf{B}$. The domain of $\widehat{\lambda}$ is a subalgebra of $\mathbf{A}$ that includes the domain of $\lambda$, which is $e(A)$, so since $e$ is dense for $\mathbf{A}$ we get that the domain of $\widehat{\lambda}$ is $\mathbf{A}$. Thus $\widehat{\lambda} \in \operatorname{Hom}(\mathbf{A}, \mathbf{B})$ and $e(\widehat{\lambda}) = \lambda$. This proves (2) (ii). □

**Theorem 2.3.** *Let $\mathcal{V}$ be a variety and let $e$ be a unary term of $\mathcal{V}$ for which $\mathcal{V} \models e^2 = e$. Let $\mathcal{S}$ denote the full subcategory of $\mathcal{V}$ whose objects are the algebras in $\mathcal{V}$ which are separated by $e$. Let $\mathcal{D}$ denote the full subcategory of $\mathcal{S}$ whose objects are the algebras in $\mathcal{S}$ for which $e$ is dense. Then*

(1) *$e(\mathcal{D}) = e(\mathcal{S})$, and this category is a prevariety, and*
(2) *$e$ is a categorical equivalence from $\mathcal{D}$ to $e(\mathcal{D})$.*

*Proof.* Since $\mathcal{D}$ is a subclass of $\mathcal{S}$, to prove the first claim of (1) we must show that for every $\mathbf{A} \in \mathcal{S}$ there is an algebra $\mathbf{A}' \in \mathcal{D}$ such that $e(\mathbf{A}') = e(\mathbf{A})$. Simply take $\mathbf{A}'$ to be the subalgebra of $\mathbf{A}$ generated by $e(A)$. Clearly this $\mathbf{A}'$ is in $\mathcal{S}$ since $\mathcal{S}$ is a prevariety, $e$ is dense for $\mathbf{A}'$ since $e(A') = e(A)$ is a generating set for $\mathbf{A}'$, and $e(\mathbf{A}') = e(\mathbf{A})$ by the properties of the $e$–construction. The fact that $e(\mathcal{S})$ is a prevariety follows from the fact that $\mathcal{S}$ is.

Now we prove (2). First we give some standard definitions associated with categorical equivalence. One says that a functor $F \colon \mathcal{C} \to \mathcal{C}'$ is **faithful** if the induced map

$$F \colon \operatorname{Hom}_{\mathcal{C}}(\mathbf{A}, \mathbf{B}) \to \operatorname{Hom}_{\mathcal{C}'}(F(\mathbf{A}), F(\mathbf{B}))$$

is injective for each pair of objects $\mathbf{A}, \mathbf{B} \in \mathcal{C}$. We say that $F$ is **full** if this induced map is surjective for all $\mathbf{A}$ and $\mathbf{B}$. We say that $F$ is **representative** if for each object $\mathbf{C}' \in \mathcal{C}'$ there is an object $\mathbf{C} \in \mathcal{C}$ such that $F(\mathbf{C}) \cong \mathbf{C}'$. We will use the following well known theorem (a proof of which can be found in [8]): A functor is a categorical equivalence if and only if it is full, faithful and representative.

Now, $e$ is a full and faithful functor from $\mathcal{D}$ to $e(\mathcal{D})$ by Lemma 2.2 (2). It is representative because the target category is defined to be the image of the functor. Thus it is a categorical equivalence. □

In fact, what the proof of Theorem 2.3 (2) shows is that if $\mathcal{K}$ is any class of similar algebras for which $e$ is a separating and dense idempotent, then $e$ is a categorical equivalence from $\mathcal{K}$ to $e(\mathcal{K})$.

**Corollary 2.4.** *Let $\mathbf{A}$ and $\mathbf{B}$ be algebras in the same language and assume that $e$ is a separating and dense idempotent for both $\mathbf{A}$ and $\mathbf{B}$. Then $\mathbf{A} \cong \mathbf{B}$ if and only if $e(\mathbf{A}) \cong e(\mathbf{B})$.*

*Proof.* The assumptions imply that $\mathbf{A}, \mathbf{B} \in \mathcal{D}$ where $\mathcal{D}$ is the subclass of the variety $\mathcal{V} = \mathcal{V}(\{\mathbf{A}, \mathbf{B}\})$ defined in Theorem 2.3. Since $e$ is a categorical equivalence from $\mathcal{D}$ to $e(\mathcal{D})$ we get that $\mathbf{A} \cong \mathbf{B}$ if and only if $e(\mathbf{A}) \cong e(\mathbf{B})$. □

## 3. Strictly Simple and Term Minimal Algebras

If $\mathbf{A}$ is an algebra, then the unary term $e$ is a **minimal idempotent** of $\mathbf{A}$ if the term operation associated to $e$ is a nonconstant idempotent operation whose range is minimal among ranges of nonconstant idempotent term operations of $\mathbf{A}$. This is a property of $e$ which can be expressed equationally, and therefore the notion can be extended from single algebras to varieties of algebras as follows. Say that $e$ is a minimal idempotent for $\mathcal{V}$ if

- $\mathcal{V} \models e^2(x) = e(x)$,
- $\mathcal{V} \not\models e(x) = e(y)$,
- If $\mathcal{V} \models f^2(x) = f(x) = ef(x)$, then $\mathcal{V} \models f(x) = e(x)$ or $\mathcal{V} \models f(x) = f(y)$.

These equations imply that $e$ interprets as a minimal idempotent or as a constant in any $\mathbf{B} \in \mathcal{V}$. If $\mathcal{V} = \mathcal{V}(\mathbf{A})$, then $e$ is a minimal idempotent for $\mathcal{V}$ if and only if it is a minimal idempotent for $\mathbf{A}$ in the sense originally defined.

We call an algebra or class of algebras **term minimal** if $e(x) = x$ is a minimal idempotent. If one localizes an algebra $\mathbf{A}$ to the neighborhood defined by a minimal idempotent one obtains a term minimal algebra, and trivially this is how all term minimal algebras arise. The class of term minimal algebras is a rich and complicated one which perhaps we will never understand. However, there are important special cases were a full description of the clones of term minimal algebras is known. One such special case is the class of algebras which are expansions of finite algebras by constants. These term minimal algebras are called *E-minimal algebras* (in [4] and [7], for instance). The classification of $E$–minimal algebras which are not strongly solvable can be found in [4], while the strongly solvable case is handled in [7]. Another understood class of term minimal algebras are the term minimal strictly simple algebras. An algebra is **strictly simple** if it is finite, simple and has no nontrivial proper subalgebras. The clones of term minimal strictly simple algebras are classified in [9]. Our first result in this section is a consequence of this classification and Corollary 2.4.

**Theorem 3.1.** *Let $\mathbf{A}$ and $\mathbf{B}$ be strictly simple algebras which generate the same variety. Then either $\mathbf{A} \cong \mathbf{B}$ or both $\mathbf{A}$ and $\mathbf{B}$ are strongly abelian.*

*Proof.* Since $\mathbf{A}$ is finite, it has a minimal idempotent; let $e$ be one such. Then $e$ is a minimal idempotent for $\mathcal{V} = \mathcal{V}(\mathbf{A}) = \mathcal{V}(\mathbf{B})$, and hence for

**B**. The local algebras $e(\mathbf{A})$ and $e(\mathbf{B})$ are term minimal, and since
$$\mathcal{V}(e(\mathbf{A})) = e(\mathcal{V}) = \mathcal{V}(e(\mathbf{B})),$$
it follows that $e(\mathbf{A})$ and $e(\mathbf{B})$ generate the same variety. By Lemma 2.1, $e$ is separating for both $\mathbf{A}$ and $\mathbf{B}$. Since $\mathbf{A}$ and $\mathbf{B}$ have no proper nontrivial subalgebras, $e$ is dense for both $\mathbf{A}$ and $\mathbf{B}$. By Corollary 2.4 we need only to prove that $e(\mathbf{A}) \cong e(\mathbf{B})$ or that both $\mathbf{A}$ and $\mathbf{B}$ are strongly abelian.

We claim that $e(\mathbf{A})$ and $e(\mathbf{B})$ are strictly simple. To see that $e(\mathbf{A})$ has no proper nontrivial subalgebras, it suffices to note that this property is inherited from $\mathbf{A}$ since (using the notation of Theorem 2.3) $\mathbf{A} \in \mathcal{D}$, $\mathcal{D}$ and $e(\mathcal{D})$ are closed under the formation of subalgebras and $e \colon \mathcal{D} \to e(\mathcal{D})$ is a categorical equivalence. If $\theta$ is any congruence on $e(\mathbf{A})$, then for $\widehat{\theta}$ equal to the transitive closure of $\mathrm{Sg}^{\mathbf{A}^2}(\theta)$ we have that $e(\widehat{\theta}) = \theta$. Since $\mathbf{A}$ is simple, we conclude that $e(\mathbf{A})$ is simple as well. This argument proves that $e(\mathbf{A})$ and $e(\mathbf{B})$ are strictly simple.

Now we can use the classification of strictly simple term minimal algebras. The crucial facts we need, which follow from the classification, are these: If $\mathbf{T}$ is a term minimal strictly simple algebra, then

(i) either $\mathbf{T}$ is the unique strictly simple algebra in $\mathcal{V}(\mathbf{T})$, or $\mathbf{T}$ is an abelian algebra with no trivial subalgebras;

(ii) if $\mathbf{T}$ is an abelian algebra, then $\mathbf{T}$ is affine or essentially unary.

Since $e(\mathbf{A})$ and $e(\mathbf{B})$ are term minimal strictly simple algebras which generate the same variety, we must have by (i) that $e(\mathbf{A}) \cong e(\mathbf{B})$ or else that $e(\mathbf{A})$ and $e(\mathbf{B})$ are both abelian and have no trivial subalgebras. The property of being essentially unary is equational, so either $e(\mathbf{A})$ and $e(\mathbf{B})$ are both affine or they are both essentially unary. It is a consequence Theorem 12.4 of [2] that two strictly simple affine algebras which generate the same variety are isomorphic. Therefore, we have that $e(\mathbf{A}) \cong e(\mathbf{B})$ or else that $e(\mathbf{A})$ and $e(\mathbf{B})$ are both essentially unary.

To finish the proof, we note that whenever $\mathbf{A}$ is a finite simple algebra and $e$ is any nonconstant polynomial of $\mathbf{A}$, then $e(A)$ contains a minimal set $N$. If $\mathbf{A}|_N$ is essentially unary, then $\mathrm{typ}\,\{\mathbf{A}\} = \{\mathbf{1}\}$ and this implies that $\mathbf{A}$ is strongly abelian. Therefore, we have that $e(\mathbf{A}) \cong e(\mathbf{B})$ or else that $\mathbf{A}$ and $\mathbf{B}$ are both strongly abelian. $\square$

One can easily find nonisomorphic strictly simple $G$–sets which generate the same variety, so the assumption in the previous corollary that $\mathbf{A}$ and $\mathbf{B}$ are not strongly abelian is necessary to show that $\mathbf{A} \cong \mathbf{B}$.

Theorem 3.1 resolves an open question from [4]. It is proved in Theorem 14.8 of [4] that if $\mathbf{A}$ and $\mathbf{B}$ are finite simple algebras which generate the same variety then they have the same tame congruence theoretic type. It is also shown that if this type is **3** or **4**, then $\mathbf{A} \cong \mathbf{B}$. Examples are given to show that if the type is **1** or **5** then it is possible that $\mathbf{A} \not\cong \mathbf{B}$. In Exercise 14.9 (3) it is asked if there exist nonisomorphic simple algebras of type **2** which generate the same variety. There do not exist such algebras, for in [10] it is shown that any simple algebra of type **2** is strictly simple. Therefore, from these remarks and the previous theorem we get that:

**Corollary 3.2.** *If $\mathbf{A}$ and $\mathbf{B}$ are simple algebras which generate the same variety and* $\mathrm{typ}\,\{\mathbf{A}\} \in \{\mathbf{2}, \mathbf{3}, \mathbf{4}\}$, *then* $\mathbf{A} \cong \mathbf{B}$.

## 4. Projective Algebras

In this section we want to use the localization functor to determine when an algebra $\mathbf{P}$ is projective. By a **category of algebras** we will mean any full subcategory of the category of all algebras in a given language. If $\mathcal{P}$ is a category of algebras, then an algebra $\mathbf{P} \in \mathcal{P}$ is **projective in** $\mathcal{P}$ if whenever

(a) $\mathbf{A}, \mathbf{B} \in \mathcal{P}$,
(b) $\sigma \colon \mathbf{A} \to \mathbf{B}$ is a surjective homomorphism and
(c) $\varphi \colon \mathbf{P} \to \mathbf{B}$ is any homomorphism,

then there exists $\bar{\varphi} \colon \mathbf{P} \to \mathbf{A}$ such that $\sigma \circ \bar{\varphi} = \varphi$. We will call $\mathbf{P}$ simply **projective** if it is projective in the variety it generates. In this section we explain how to simplify the test for whether an algebra $\mathbf{P}$ is projective in a given variety of algebras.

Let $\mathcal{K}$ be a class of similar algebras and assume that $\mathbf{P}$ is in the variety generated by $\mathcal{K}$. It is well known and easy to prove that the following conditions are equivalent:

(1) $\mathbf{P}$ is projective in $\mathcal{V}(\mathcal{K})$;
(2) $\mathbf{P}$ is a retract of a free algebra in $\mathcal{V}(\mathcal{K})$;
(3) $\mathbf{P}$ is a retract of a free algebra in $\mathsf{ISP}(\mathcal{K})$;
(4) $\mathbf{P}$ is projective in $\mathsf{ISP}(\mathcal{K})$.

We want to connect these conditions with projectivity in a category smaller than $\mathsf{ISP}(\mathcal{K})$. Choose and fix an $e$ which is a separating idempotent for all algebras in $\mathcal{K}$. For this fixed $e$ let $\mathcal{E}$ denote the full subcategory of $\mathsf{ISP}(\mathcal{K})$ whose object class consists of the algebras in $\mathsf{ISP}(\mathcal{K})$ for which $e$ is dense. Call $\mathbf{P}$ a **dense projective** in $\mathcal{V}(\mathcal{K})$ if $\mathbf{P}$ is projective in $\mathcal{V}(\mathcal{K})$ and $e$ is dense for $\mathbf{P}$.

**Theorem 4.1.** *The dense projectives in $\mathcal{V}(\mathcal{K})$ are precisely the algebras which are projective in $\mathcal{E}$.*

*Proof.* Any algebra projective in $\mathcal{V}(\mathcal{K})$ lies in $\mathsf{ISP}(\mathcal{K})$, as we have observed, and therefore any dense projective must lie in $\mathcal{E}$. If $\mathbf{P}$ is such an algebra (projective in $\mathcal{V}(\mathcal{K})$ and lying in $\mathcal{E}$), then $\mathbf{P}$ is projective in $\mathcal{E}$ since $\mathcal{E}$ is a full subcategory of $\mathcal{V}(\mathcal{K})$. Therefore, what we must prove is that $\mathcal{E}$ has no projectives other than these. That is, we must show that if $\mathbf{P} \in \mathcal{E}$ is projective in $\mathcal{E}$, then it is projective in $\mathcal{V}(\mathcal{K})$. To prove this, we will show that any $\mathbf{P}$ which is projective in $\mathcal{E}$ is a retract of a free algebra of $\mathcal{V}(\mathcal{K})$.

Let $\mathbf{F}$ be a free algebra of rank large enough for there to be a surjection $\sigma\colon \mathbf{F} \to \mathbf{P}$. Let $\mathbf{L}$ be the subalgebra of $\mathbf{F}$ generated by $e(F)$. Since $\mathbf{F} \in \mathsf{ISP}(\mathcal{K})$, we get that $\mathbf{L} \in \mathsf{ISP}(\mathcal{K})$ and of course $\mathbf{L}$ is generated by $e(F) = e(L)$. Hence $\mathbf{L} \in \mathcal{E}$. Furthermore, since $\mathbf{L}$ is generated by $e(L)$, the homomorphism $\sigma|_L\colon \mathbf{L} \to \mathbf{P}$ is surjective: indeed, the image of $\sigma|_L$ is

$$\begin{aligned}\sigma|_L(\mathbf{L}) &= \sigma(\mathrm{Sg}^{\mathbf{F}}(e(F)))\\ &= \mathrm{Sg}^{\mathbf{P}}(\sigma(e(F)))\\ &= \mathrm{Sg}^{\mathbf{P}}(e(\sigma(F)))\\ &= \mathrm{Sg}^{\mathbf{P}}(e(P))\\ &= \mathbf{P}\end{aligned}$$

because $\mathbf{P}$ is in $\mathcal{E}$. Since $\mathbf{P}$ is projective in $\mathcal{E}$, there is a homomorphism $\tau\colon \mathbf{P} \to \mathbf{L}$ such that $\sigma|_L \circ \tau = \mathrm{id}_{\mathbf{P}}$. If $\iota\colon \mathbf{L} \to \mathbf{F}$ denotes inclusion, then $\iota \circ \tau\colon \mathbf{P} \to \mathbf{F}$ is a homomorphism for which we have

$$\sigma \circ (\iota \circ \tau) = (\sigma \circ \iota) \circ \tau = \sigma|_L \circ \tau = \mathrm{id}_{\mathbf{P}}.$$

Hence $\iota \circ \tau$ is a right inverse for $\sigma$ which shows that $\mathbf{P}$ is a retract of $\mathbf{F}$. This finishes the proof of the theorem. □

**Theorem 4.2.** *Assume that $e$ is a separating idempotent for all algebras in $\mathcal{K}$ and that $e$ is separating and dense for some $\mathbf{P} \in \mathcal{V}(\mathcal{K})$. Then $\mathbf{P}$ is projective in $\mathcal{V}(\mathcal{K})$ if and only if $e(\mathbf{P})$ is projective in $\mathcal{V}(e(\mathcal{K}))$.*

*Proof.* Let $\mathcal{E}$ be the full subcategory of $\mathsf{ISP}(\mathcal{K} \cup \{\mathbf{P}\})$ whose object class consists of the algebras for which $e$ is dense. Under the given hypotheses we have that $\mathbf{P} \in \mathcal{E}$, so by Theorem 4.1 $\mathbf{P}$ is projective in $\mathcal{V}(\mathcal{K})$ if and only if it is projective in $\mathcal{E}$. By Theorem 2.3 and the remarks immediately after its proof, the localization functor is a categorical equivalence from $\mathcal{E}$ to $e(\mathcal{E})$. Therefore, the following are equivalent:

(1) $\mathbf{P}$ is projective in $\mathcal{V}(\mathcal{K})$;
(2) $\mathbf{P}$ is projective in $\mathcal{E}$;
(3) $e(\mathbf{P})$ is projective in $e(\mathcal{E}) = \mathsf{ISP}(e(\mathcal{K}) \cup \{e(\mathbf{P})\}) \supseteq \mathsf{ISP}(e(\mathcal{K}))$;
(4) $e(\mathbf{P})$ is projective in $\mathcal{V}(e(\mathcal{K}))$.

This finishes the proof. □

We apply Theorem 4.2 to describe the projective strictly simple algebras. In this case, we take $\mathcal{K} = \{\mathbf{A}\}$ where $\mathbf{A}$ is a strictly simple algebra. Let $e$ be any minimal idempotent for $\mathbf{A}$. As we observed in the proof of Theorem 3.1, $e$ is a separating and dense idempotent for $\mathbf{A}$ and $e(\mathbf{A})$ is strictly simple. By Theorem 4.2 we have that $\mathbf{A}$ is projective if and only if $e(\mathbf{A})$ is projective. Therefore the characterization of projective strictly simple algebras reduces to the term minimal case, which is handled in the next lemma.

**Lemma 4.3.** *A strictly simple term minimal algebra is projective if and only if it is not definitionally equivalent to an irregular $G$–set.*

*Proof.* Most of this is already proved in Corollary 2.7 of [6]. There it is shown that if $\mathbf{T}$ is nonabelian or has a trivial subalgebra then $\mathbf{T}$ is projective. The remaining cases to consider are when $\mathbf{T}$ is abelian and has no trivial subalgebras. Moreover, the proof of Corollary 2.7 in the subcase where every element of $\mathbf{T}$ is the interpretation of a constant term works equally well in the abelian and nonabelian cases. Thus, the only strictly simple term minimal algebras which could fail to be projective are the abelian ones which have no proper subalgebras and no constant terms. As we have mentioned previously, the abelian term minimal strictly simple algebras are essentially unary or affine. If $\mathbf{T}$ is essentially unary, then it is equivalent to a (transitive) primitive $G$–set. If $\mathbf{T}$ is affine, then it is the expansion by translations of a simple affine module.

In the affine case $\mathbf{T}$ is projective. In fact, in this case the collection of unary term operations of $\mathbf{T}$ coincides with the group of additive translations $\{x \mapsto x + t \mid t \in T\}$. The automorphism group of $\mathbf{T}$ also coincides with this group. The coincidence of these groups implies that $\mathbf{T} \cong \mathbf{F}_{\mathcal{V}(\mathbf{T})}(1)$. Since $\mathbf{T}$ is free it is projective.

In the unary case, assume that $\mathbf{T}$ is a faithful primitive $G$–set. Since $\mathbf{T}$ is faithful, the left regular representation of $G$, call it $\mathbf{L}$, belongs to $\mathcal{V}(\mathbf{T})$. Furthermore, since $\mathbf{T}$ is transitive and $\mathbf{L}$ is free there is a surjective homomorphism from $\mathbf{L}$ onto $\mathbf{T}$. If $\mathbf{T}$ is projective, then it must be that $\mathbf{T}$ is a retract of $\mathbf{L}$. But $\mathbf{L}$ has no proper subalgebras so we must have that $\mathbf{T} \cong \mathbf{L}$. Thus, if $\mathbf{T}$ is projective then it must be equivalent to a regular $G$–set. Conversely, if $\mathbf{T}$ is equivalent to a regular $G$–set then, since $\mathbf{T}$ has no proper subalgebras, $\mathbf{T}$ must be isomorphic to $\mathbf{L}$. Since $\mathbf{L} \cong \mathbf{T}$ is free, we get that $\mathbf{T}$ is projective. $\square$

**Corollary 4.4.** *A strictly simple algebra $\mathbf{A}$ is projective if and only if it has no idempotent unary term $e$ such that $e(\mathbf{A})$ is definitionally equivalent to an an irregular $G$–set.*

This corollary can be used to give a new proof of Theorem 3.1, if one uses the easily proven fact that two subdirectly irreducible projective algebras which generate the same variety must embed into one another.

We actually have a good deal more information than we have stated about the class of projectives in $\mathcal{V}(\mathbf{A})$ when $\mathbf{A}$ is a strictly simple algebra. If $e$ is a fixed minimal idempotent of $\mathbf{A}$, then an algebra $\mathbf{P} \in \mathcal{V}(\mathbf{A})$ is a dense projective (recalling terminology from Theorem 4.1) if and only if $e$ is separating and dense for $\mathbf{P}$ and $e(\mathbf{P})$ is projective in $\mathcal{V}(e(\mathbf{A}))$. Therefore the determination of dense projectives in $\mathcal{V}(\mathcal{K})$ can be reduced to the determination of the projectives in $\mathcal{V}(e(\mathcal{K}))$. The entire class of projectives in $\mathcal{V}(e(\mathbf{A}))$ is not hard to describe when $\mathbf{A}$ is abelian. When $\mathbf{A}$ is nonabelian, the class of projectives in $\mathcal{V}(e(\mathbf{A}))$ has been worked out in some key cases. For example, the full class of projectives for the varieties of Boolean algebras, distributive lattices and semilattices can be found respectively in [3], [1] and [5]. It seems a (difficult but) feasible project to characterize the dense projectives in varieties generated by strictly simple algebras using Theorem 4.2.

## References

[1] R. Balbes and A. Horn, *Projective distributive lattices*, Pac. J. Math. **33** (1970), 273–279.

[2] R. Freese and R. McKenzie, *Commutator Theory for Congruence Modular Varieties*, LMS Lecture Notes v. 125, Cambridge University Press, 1987.

[3] R. Haydon, *On a problem of Pełczyński, Milutin spaces, Dugundji spaces and AE (0–dim)*, Studia Math. **52** (1974), 23–31.

[4] D. Hobby and R. McKenzie, *The Structure of Finite Algebras*, Contemporary Mathematics v. 76, American Mathematical Society, 1988.

[5] A. Horn and N. Kimura, *The category of semilattices*, Algebra Universalis **1** (1971), 26–38.

[6] K. A. Kearnes and Á. Szendrei, *A characterization of minimal locally finite varieties*, to appear in the Transactions of the AMS.

[7] E. W. Kiss, *An easy way to minimal algebras*, to appear in the International Journal of Algebra and Computation.

[8] S. Mac Lane, *Categories for the Working Mathematician*, Graduate Texts in Mathematics vol. 5, Springer–Verlag, 1971.

[9] Á. Szendrei, *Term minimal algebras*, Algebra Universalis **32** (1994), 439–477.

[10] M. Valeriote, *Finite simple abelian algebras are strictly simple*, Proceedings of the AMS **108** (1990), 49–57.

(Keith A. Kearnes) DEPARTMENT OF MATHEMATICAL SCIENCES, UNIVERSITY OF ARKANSAS, FAYETTEVILLE, AR 72701, USA.

(Ágnes Szendrei) BOLYAI INSTITUTE, ARADI VÉRTANÚK TERE 1, H–6720 SZEGED, HUNGARY.
  *E-mail address*: `kearnes@comp.uark.edu`
  *E-mail address*: `a.szendrei@math.u-szeged.hu`